\documentstyle[12pt]{article}
\normalbaselines
\catcode `\@=11
\@addtoreset{equation}{section}

\def\section{\@startsection {section}{1}{\z@}{-3.5ex plus -1ex minus
     -.2ex}{2.3ex plus .2ex}{\normalsize\bf}}
\def\subsection{\@startsection{subsection}{2}{\z@}{-3.25ex plus -1ex minus
 -.2ex}{1.5ex plus .2ex}{\normalsize\bf}}
\def\thebibliography#1{\section*{References\markboth
  {REFERENCES}{REFERENCES}}\list
  {[\arabic{enumi}]}{\settowidth\labelwidth{[#1]}\leftmargin\labelwidth
  \advance\leftmargin\labelsep
  \usecounter{enumi}}
  \def\newblock{\hskip .11em plus .33em minus -.07em}
  \sloppy
  \sfcode`\.=1000\relax}
 
\catcode `\@=12
\begin{document}
\vspace*{2.5cm}
\begin{flushright} 
SNUTP 97-056 
\end{flushright}
\begin{center}
{ \bf The $q$-Deformed Oscillator Representations and }\\
\vspace{0.5cm}
{\bf  Their Coherent States of the $su(1,1)$ Algebra}
\vspace{1.0cm}\\
\end{center}
\hspace*{1in}
\begin{minipage}{13cm}
Phillial Oh$^{1}$ and Chaiho Rim$^{2}$ \vspace{0.3cm}\\
 $^{1}$ Department of Physics, Sung Kyun Kwan University\\
\makebox[3mm]{ }Suwon, 440-746,  Korea \\
\makebox[3mm]{ }E-mail address: ploh@newton.skku.ac.kr\\
 $^{2}$ Department of Physics, Chonbuk National University\\
\makebox[3mm]{ }Chonju, 561-756,  Korea\\
\makebox[3mm]{ }E-mail address: rim@phy0.chonbuk.ac.kr
\end{minipage}
 
\vspace*{0.5cm}
\begin{abstract}
\noindent
We present various oscillator representations of the 
$q$-deformed $su(1,1)$ algebra such as  the
Holstein-Primakoff, the Dyson, the Fock-Bargmann, the anyonic, 
and the parabose oscillator representations
and discuss their coherent states with the resolution of unity. 
\end{abstract}

\section{\hspace{-4mm}.\hspace{2mm} Introduction}

The purpose of this work is to
study the various $q$-deformed oscillator 
representations of the $su_q(1,1)$
algebra \cite{jimbo,bied_mac} each of which has its own merit.
It is rather well known that 
Holstein-Primakoff (HP) and Dyson  representation 
in the $su(2)$ case have been useful, among others,
in studying the magnon \cite{kitt}
and the nuclear excitation \cite{rowe}. 
On the other hand, the Fock-Bargmann (FB) representation provides
a convenient tool in describing
the quantum system in terms of path-integral formalism. 
In addition, the parabose oscillator system \cite{green}
plays an important role in studying 
the exclusion statistics \cite{haldane} of  quasi-particles  
through the Calogero model \cite{calo}.

In this paper, we construct various oscillator  representations
and their coherent states
of the $q$-deformed $su(1,1)$ algebra and  discuss
the relation among them. 
In section 2, we present the $q$-deformation of
the HP, the Dyson, the anyonic, and the 
FB representations
starting from the coadjoint orbit approach.
The oscillator representations and their 
coherent states  are given.
In section 3,  the quadratic oscillator system with the 
Casimir invariant $-{3 \over 16}$ is considered. 
 We give its generalization, the parabose oscillator 
system, and the coherent state. Section 4 is the conclusion.

\section{\hspace{-4mm}.\hspace{2mm} 
	The HP and  Other Representations }

\subsection{\hspace{-5mm}.\hspace{2mm}
The Coadjoint Orbit Approach }

We start with a brief description of the
coadjoint orbit of the group $SU(1,1)$ because
it gives a clear idea about the 
symplectic structure and canonical
quantization of the phase space of HP representation.
Let us express the group element $g$ of $SU(1,1)$ which is 
a $2 \times 2$ matrix by the kets $(\vert Z_1>, \vert Z_2>)$ 
with $ \vert Z_p>=(Z_p^{(1)},Z_p^{(2)})^T$.
We can get rid of $\vert Z_2>$ by making use 
of the special unitary property of $SU(1,1)$
and express $g$ in terms of  $ \vert Z_1>\equiv \vert Z> $ only.
The isospin charge on the coadjoint orbit $SU(1,1)/U(1)$ 
becomes  \cite{oh955}
$\tilde K_a=-2\mbox{Tr}(gxg^{-1}T_a)
=-2ik_0 < Z \vert T_a \vert Z>$
where $x=ik_0~\mbox{diag}(-1,1)$ 
and $k_0$ is restricted to be a positive integer or half odd integer. 
The $T_a$'s are generators of the  $su(1,1)$ algebra.
Employing the projective coordinates, 
$\vert Z > =\frac{1}{1-\vert\xi\vert^2} (1, \xi)^T $,
we have  the  isospin charge  given by 
\begin{equation}
\tilde K_1=k_0\frac{\xi+\bar\xi}{1-\vert\xi\vert^2},~~
\tilde K_2=ik_0\frac{\xi-\bar\xi}{1-\vert\xi\vert^2},~~
\tilde K_3=k_0\frac{1+\vert\xi\vert^2}{1-\vert\xi\vert^2}\,.
\label{isospinclass}
\end{equation}

The symplectic structure of the phase space can be used to
construct the action
\begin{equation}
L=2\mbox{Tr}(xg^{-1}\dot g)-H(\tilde K_a)
=i k_0\frac{\dot{\bar \xi} \xi-{\bar \xi}\bar\xi}
{1-\vert\xi\vert^2}-H(\tilde K_a)\,.
\label{action}
\end{equation}
On this phase space one can demonstrate that the isospin 
charges satisfy the $su(1,1)$ Lie-Poisson algebra, 
$\{\tilde K_a,  \tilde K_b\}=-\epsilon_{ab}^{\ \ c} \tilde K_c$.
(We are using the metric $\eta^{ab}=(1,1,-1)$ in lowering and 
raising the index.)
Transforming the action Eq. (\ref{action}) into a canonical form 
by the use of the Darboux variables \cite{oh945}, we have 
$L_\alpha=\frac{i}{2}({\bar \alpha}\dot \alpha-\dot{\bar \alpha}
\alpha)-H(\tilde K_a)$, 
and we have the canonical 
the Poisson bracket for $\alpha$,  $\{\alpha,\bar\alpha\}=i$.
The  isospin charges (\ref{isospinclass}) can be expressed as follows:
\begin{equation}
\tilde K_1=\frac{(\alpha+\bar\alpha)}
{2}\sqrt{2k_0+\vert\alpha\vert^2},~
\tilde K_2=\frac{i(\alpha-\bar\alpha)}
{2}\sqrt{2k_0+\vert\alpha\vert^2},~
\tilde K_3=\vert\alpha\vert^2+k_0.
\end{equation}

\subsection{\hspace{-5mm}.\hspace{2mm}
The Oscillator Representation and the Coherent State }

The HP realization in terms of quantum mechanical oscillators 
is obtained if we replace the Poisson bracket for $\alpha$
with the Dirac commutator $ [a_- , a_+ ]=1$,
$\alpha \rightarrow a_- ,~ \bar\alpha \rightarrow a_+$, 
and perform the normal ordering of the resulting
operators:
\begin{equation}
K_+ = a_+ \sqrt{2k_0+a_+ a_-},~
K_- =\sqrt{2k_0+a_+  a_-}a_-,~
K_0= a_+ a_- + k_0.
\end{equation}
If we shift the square root part in $K_-$ into the $K_+$ side, 
we get the Dyson realization:
\begin{equation}
K_+ =a_+  (2k_0+a_+ a_- ),~
K_- =a_-,~
K_0 = a_+ a_- +k_0\,.
\end{equation}
If we use the holomorphic coordinate, 
$a_- \rightarrow d/d\xi,~~ a_+ \rightarrow \xi$, 
then we get the FB realizations:
\begin{equation}
K_+ (\xi) =  \xi^2 {d \over d \xi} + 2k_0  \xi\,~~
K_-(\xi) = {d \over d \xi},~~
K_0 (\xi) = \xi {d \over d \xi} +k_0.
\label{fbfb}
\end{equation}
The generators satisfy the $su(1,1)$ algebra
\begin{equation}
[ K_0,  K_{\pm}]=\pm K_{\pm}, \quad 
[ K_+,  K_-]=-2K_0
\end{equation}
and the Casimir invariant is expressed as
$ C=K_0(K_0-1)-K_+K_-\,$.

We define the coherent state of the $SU(1,1)$ as \cite{klau}
\begin{equation}
\vert \xi>=\exp(\xi  K_+)\vert \Lambda>,~~ h\vert \Lambda>=
e^{i\psi(h)}\vert\Lambda>,~~ h\in U(1).
\end{equation}
Here $\vert\Lambda>$ is the lowest weight vector (the vacuum).
Once we require 
the $K_-$ to be conjugate to the $K_+$, 
then  the resolution of unity is given 
in terms of the Liouville measure
independent of  the representations we considered above:
\begin{equation}
I=\int d\bar\xi d\xi \vert\xi>
<\bar\xi\vert  G(|\xi|),~~
G(|\xi|)=(1-\vert\xi\vert^2)^{2k_0 -2}.
\label{unity}
\end{equation}
On the other hand,
one can also define the coherent state of the oscillator algebra,
the Glauber coherent state, in the HP case, because
the $a_-$ and the $a_+ $ are manifestly conjugate to each other: 
\begin{equation}
|\xi] = e ^{\xi a_+ } |0>\,.
\end{equation}
The resolution of  unity is satisfied with the measure of
Eq. (\ref{unity}) replaced by
$G(|\xi|) = {1 \over \pi} e^{- |\xi|^2}$.

\subsection{\hspace{-5mm}.\hspace{2mm}
The $q$-Deformation of $su(1,1)$}

The $q$-deformed algebra $su_q(1,1)$ is given as \cite{jimbo}
\begin{equation}
[Q_0, Q_{\pm}]=\pm  Q_{\pm}\,,\quad
[ Q_+,  Q_-]=-[2Q_0]_{q}\,
\end{equation}
where the $q$-deformation is defined as
$[x]_q\equiv \frac{q^x-q^{-x}}{q-q^{-1}}$.
The $q$-deformed Casimir invariant is given by 
$C_q = [Q_0]_q[Q_0 - 1]_q - Q_+ Q_-$.
One can obtain the explicit form of the $q$-deformed 
generators following  Ref. \cite{curt}:
\begin{equation}
Q_0=K_0\,,\quad 
Q_-=K_-f(K_0)\,,\quad
Q_+ =f(K_0)K_+\,
\label{q11_rz}
\end{equation}
with 
$f(K_0)= \sqrt{[K_0-k_0]_q [K_0 + k_0 -1]_q
\over (K_0 -k_0) (K_0 + k_0 -1)}$.
We used the symmetric form of the $Q_{\pm}$ 
to make them manifestly conjugate to each other.
However, this is not the unique choice (see below).

\subsection{\hspace{-5mm}.\hspace{2mm}
The Biedenharn Oscillator and the Coherent State }

According to Eq. (\ref{q11_rz}), 
one can assume that the generators of the $su_q(1,1)$ algebra 
can be expressed by 
\begin{equation}
Q_0= K_0 \equiv N+ k_0\,,~~
Q_- = (a_q)_- \sqrt{[N +2k_0 -1]_q} \,,~~
Q_+ = \sqrt{[N + 2k_0 -1]_q} (a_q)_+\,
\end{equation}
where $(a_q)_\pm$ is identified as 
\begin{equation}
(a_q)_- = a_-  \sqrt{[N]_q \over N}\,,\quad
(a_q)_+ =  \sqrt{[N]_q \over N}\,\, a_+ \,.
\end{equation}
The $(a_q)_\pm$ satisfy  the $q$-deformed 
oscillator algebra of the Biedenharn type
\begin{equation}
(a_q)_- (a_q)_+ - q (a_q)_+ (a_q)_- = q^{-N}\,.
\label{Bied_com}
\end{equation}

Let us define the coherent state of the 
$q$-deformed  $su(1,1)$  by 
\begin{equation}
|z>  = e_q^{\bar z Q_+} |0>
\label{q11_coherent}
\end{equation}
where the   $q$-deformed exponential function is given by
$e_q^x \equiv \sum_{n=0}^{\infty} {x^n \over [n]_q!}$.
The resolution of unity is given as
\begin{equation}
I = \int d^2_q z \, g(z)\, |z><z|\,,
\end{equation}
and the measure  $g(z)$ is given as \cite{oh_rim}
\begin{equation}
g(|z|)= \left \{ \begin{array}{ll}
{[2k_0 -1]_q \over \pi} (1 - |z|^2)_q^{2k_0 -2}&
 \mbox{for } ~k_0= \mbox{integer} >1 \\
 {|z|^2 \over \pi} & \mbox{for }~ k_0= 1
\end{array}\right.
\label{q11_measure}
\end{equation}
where the $q$-deformed function is defined as 
$(1-x)^n _q = \sum_{m=0}^n {[n]_q! \over [m]_q! [n-m]_q!} (-x)^m$.
(For $k_0= {1 \over 2}$, see below Eq. (\ref{aq_measure})).
Here,  the two dimensional integration is defined as 
$d^2_q z \equiv {1 \over 2} d\theta\,\, d_q |z|^2 $.
The angular integration is an ordinary integration, 
$0 \le \theta \le 2\pi$.
The radial part is a  $q$-integration which is the inverse 
operation of $q$-derivative defined as
\begin{equation}
{d \over d_q z} f(z) = {f(qz) - f(q^{-1}z) \over z(q- q^{-1})}\,.
\label{q_derivative}
\end{equation}
One can check that $I$ commutes with the $su_q(1,1)$ generators.

As in the case with no deformation,  
one can define the coherent state of 
the $q$-deformed oscillator by
$|z] = e_q^{z a_+} |0>$, 
and the measure for the resolution of unity is given as 
\begin{equation}
g(| z|) = {1 \over \pi} e_q ^{-|z|^2}\,.
\label{aq_measure}
\end{equation} 
Note that the domain is over an infinite plane unlike in the 
$su(1,1)$ coherent state and also 
this holds for any value of $k_0$.
Therefore,  for $k_0={1 \over 2}$, 
the $q$-analogue of Glauber coherent state is suitable. 
For other value of $k_0$, 
one can use the Bargmann measure defined in 
Eq.~(\ref{aq_measure}) or  the Liouville measure in 
Eq.~(\ref{q11_measure})  depending
on the choice of the coherent state.

\subsection{\hspace{-5mm}.\hspace{2mm}
The  Macfarlane Oscillator and the Coherent State 
 }
Let  us  assume that the  $su_q(1,1)$ algebra is given by 
\begin{equation}
Q_0= N+ k_0\,,~~
Q_- = (b_q)_- \sqrt{q^{-(N-1)} [N +2k_0 -1]_q} \,,~~
Q_+=\sqrt{q^{-(N-1)} [N + 2k_0 -1]_q} (b_q)_+\,
\end{equation}
where the  $b_q$'s are the Macfarlane oscillators: 
\begin{equation}
(b_q)_- = (a_q)_- q^{N-1 \over 2}
	= a_- \sqrt{\{N\}_q \over N}\,,\quad
(b_q)_+ = q^{N-1 \over 2}(a_q)_+
	= \sqrt{\{N\}_q \over N} a_+ \,.
\end{equation}
In the above we introduced a new definition of $q$-number,
\begin{equation}
\{x\}_q = {q^{2x} - 1 \over q^2 -1} = [x]_q \,\, q^{x-1}\,.
\end{equation}
This oscillator realization gives the $q$-deformed 
oscillator algebra of Macfarlane type:
\begin{equation}
(b_q)_- (b_q)_+ - q^2 (b_q)_+ (b_q)_- = 1\,.
\label{Mac_com}
\end{equation}

If  we define the coherent state 
for $su_q(1,1)$ as  in Eq.~(\ref{q11_coherent}),
then  the same measure $g(z)$ 
in Eq.~(\ref{q11_measure})
is used for the resolution of unity. 
However, if we define another version of $q$-deformed 
Glauber coherent state as
$|z\} = E_q^{\bar z (b_q)_+}|0> $
where $E_q^x$ differs from $e_q^x$ in that 
the $q$-number $[n]_q$ is replaced by $\{n\}_q$,
$E_q^x = \sum_{n=0}^\infty {x^n \over \{n\}_q!}$, 
then the measure for the resolution of unity is given as
\begin{equation}
g(|z|) = {1 \over \pi} E_q^{- |z|^2}\,,
\label{q_Mmeas}
\end{equation}
and the domain is over an infinite plane. 

\subsection{\hspace{-5mm}.\hspace{2mm}
The $q$-Anyonic Oscillator Realization}

Let us assume  $Q_0 = K_0$ and  $Q_{\pm}$ as 
\begin{equation}
Q_- = (A_q)_- \sqrt{ (A_q)_+ (A_q)_+ +
 2[k_0 - {1 \over 2}]_q} \,,~~
Q_+ = \sqrt{ (A_q)_+ (A_q)_+ 
+ 2[k_0 - {1 \over 2}]_q}\,\, (A_q)_+ \,
\end{equation}
where 
$(A_q)_-= a_- \,\, \sqrt{[N+ k_0 - {1\over 2} ]_q - 
   [k_0 - {1 \over 2}]_q \over N}$
and $(A_q)_+$ is the  conjugate of  $(A_q)_-$.
Its commutation relation looks complicated,
\begin{equation}
((A_q)_- (A_q)_+ + [k_0 - {1 \over 2}]_q)
-q ((A_q)_+ (A_q)_- + [k_0 - {1 \over 2}]_q)
= q^{-(N + k_0 + { 1\over 2})}\,.
\end{equation}
However, the meaning of this commutation relation 
becomes clear if 
we rewrite the relation in the Macfarlane's form
\begin{equation}
(B_q)_- (B_q)_+ - q^2 (B_q)_+ (B_q)_- =1 \,,
\end{equation}
by identifying 
\begin{eqnarray}
(B_q)_- (B_q)_+ = q^{N + k_0 - {1 \over 2}}
		\,\,((A_q)_- (A_q)_+ + [k_0 - {1 \over 2}]_q)\,,
\nonumber\\
(B_q)_+ (B_q)_- = q^{N + k_0 - {3 \over 2}}
		\,\,((A_q)_+ (A_q)_- + [k_0 - {1 \over 2}]_q)\,.
\label{Bquad}
\end{eqnarray}
The conjugate relation between $(B_q)_-$ and $(B_q)_+$
can be seen  formally at the operator level in Eq.~(\ref{Bquad})
since the conjugate relation between $(A_q)_-$ and $(A_q)_+$ does
hold. However, the vacuum $|0>$ is not annihilated by $(B_q)_-$ 
unless $k_0 = {1 \over2}$, since 
\begin{equation}
(B_q)_+ (B_q)_- |0> = \{k_0 - {1 \over 2}\}_q |0>\,.
\end{equation}

This feature reflects the fact that this realization corresponds 
to the $q$-deformed non-trivial one dimensional analogue of 
the anyon which appears in two dimensional oscillator representation
with $k_0- {1 \over 2}$ being related with the statistical
parameter in the anyon physics \cite{chorim}. 
However, the fact that $(B_q)_-$ does not annihilate the 
vacuum $|0>$ implies that one cannot define a proper Hilbert space
unless extra truncation  mechanism 
of the state below the vacuum  is imposed.
It is not known if the measure for the Glauber coherent state 
of the $B_q$ oscillator system can be defined.

\subsection{\hspace{-5mm}.\hspace{2mm}
The $q$-Deformed FB Realization }

Let us consider a realization,
\begin{equation}
Q_0=K_0 = N+ k_0\,,\quad 
Q_-=K_- {[N]_q \over N} \,,\quad
Q_+= {[N + 2k_0 -1]_q \over (N+ 2k_0 -1)} K_+\,.
\label{q11_FB}
\end{equation}
These generators act on a ket  $|n>$ as 
$Q_- |n> =  [n]_q |n-1>$ 
and  $Q_+ |n> = [n+ 2k_0]_q\, |n+1>$.
Using  $<\xi|n>= \xi^n$, we have a holomorphic realization:
\begin{equation}
\hat Q_+ (\xi) =  \xi [\xi {d \over d \xi} + 2 k_0 ]_q \,,\quad
\hat Q_-(\xi) = {d \over d_q \xi} \,,\quad
\hat Q_0 (\xi) = \xi {d \over d \xi} + k_0 \,.
\label{homo}
\end{equation}
The $q$-derivative in $\hat Q_-$ is 
defined in Eq. (\ref{q_derivative}).
This implies that the oscillator realization is given as
\begin{equation}
(a_q)_-(\xi) = {d \over d_q \xi} \,,\quad
(a_q)_+(\xi) = \xi \,
\label{axi}
\end{equation}
which satisfies the $q$-deformed oscillator algebra 
of the Biedenharn type as in Eq.~(\ref{Bied_com}).

We may consider  a slightly  modified 
version of Eq.~(\ref{q11_FB}):
\begin{equation}
Q_0=K_0 = N+ k_0\,,\quad 
Q_-=K_- {[N]_q \over N}q^{N-1} \,,\quad
Q_+= q^{-(N-1)}{[N + 2k_0 -1]_q \over (N+ 2k_0 -1)} K_+\,.
\label{q11_FBM}
\end{equation}
Then the generators act on the ket as 
\begin{equation}
Q_- |n> =  \{n\}_q |n-1> \,,\qquad
Q_+ |n> = q^{-n} [n+ 2k_0]_q\, |n+1>\,.
\end{equation}
We have a similar  
$q$-deformation of the FB representation
of the $su(1,1)$ algebra as in Eq. (\ref{homo}):
\begin{equation}
\hat Q_0 (\xi) = \xi {d \over d \xi} + k_0 \,,\quad
\hat Q_-(\xi) = {D \over D_q \xi} \,,\quad
\hat Q_+ (\xi) =  \xi q^{-(2 \xi {d \over d \xi} + 2k_0 -1)}
	\{\xi {d \over d \xi} + 2 k_0 \}_q \,.
\end{equation}
The derivative in $\hat Q_-$ is replaced by 
a new  $q$-derivative
which is given by 
\begin{equation}
{D \over D_q z} f(z) = {f(q^2z) - f(z) \over z(q^2- 1)}\,.
\end{equation}
This implies that the oscillator realization is given by
\begin{equation}
(b_q)_-(\xi) = {D \over D_q \xi} \,,\quad
(b_q)_+(\xi) = \xi \,
\end{equation}
which satisfies the $q$-deformed oscillator algebra 
of the Macfarlane type as in  Eq.~(\ref{Mac_com}).
These Fock-Bargmann realizations will be useful in path-integral 
formalism.

\section{\hspace{-4mm}.\hspace{2mm} 
The Quadratic Oscillator Representation}
  
\subsection{\hspace{-5mm}.\hspace{2mm}
The Biedenharn Type and the Parabose Oscillator}

The generators of  the $su(1,1)$ algebra  is parity even  when 
they are represented in terms of the quadratic oscillators and  
$k_0$ in Eq.~(\ref{fbfb}) is $\frac{1}{4}$ or $\frac{3}{4}$.
Its $q$ deformation is 
obtained by $su_{q^2}(1,1)$ rather than $su_q(1,1)$.
The oscillator realization is obtained  in terms of  the 
Biedenharn oscillator:
\begin{equation}
Q_0 = \frac{1}{2} (N + {1 \over 2})\,,\quad
Q_- = \frac{1}{[2]_q} (a_q)_-^2 \,,\quad
Q_+ = \frac{1}{[2]_q} (a_q)_+^2 \,.
\label{Q}
\end{equation}

This quadratic oscillator representation allows the
generalized form of the algebra if the 
oscillator algebra is generalized to the  para-bose oscillator 
algebra whose explicit coordinate representation 
results in the  Calogero model 
\cite{poly_brink,cho}:
\begin{equation}
[A_-, A_+] = 1 + 2l M.
\end{equation}
Here, $l> -\frac{1}{2}$ and $M$ is the parity operator which
anti-commutes with $A_-$ and $A_+$. 
Its $q$-deformation is obtained 
if we replace the $a_-$ and the $a_+$ 
of the $Q_{\pm}$ and the $Q_0$  by  the $A_-$ and the $A_+$
in Eq. (\ref{Q}). For $f(k_0)$ in Eq.~(\ref{q11_rz}),  
$k_0= {3 \over 4} - {1 \over 2}l M$.
The Casimir invariant is given by 
$C_{q^2} = [ {1 \over 4} + {l M \over 2}]_{q^2}
[ -{3 \over 4} + {l M \over 2}]_{q^2}$
and the $q$-deformed para-bose oscillator can be expressed as
\begin{equation}
(A_q)_- = A_- \sqrt{{[A_+ A_-]_q \over A_+ A_-}},
\quad
(A_q)_+ 
=\sqrt{{[A_+ A_-]_q \over A_+ A_-}}A_+.
\end{equation}

\subsection{\hspace{-5mm}.\hspace{2mm}
The Coherent State of the Parabose Oscillator}

For the coherent state of the $su(1,1)$, we introduce the 
Fock space with the vacuum $|0>$ which is parity even 
and annihilated 
by the $K_-$. Since there are two classes of states,  
parity even and odd,
we denote them as 
\begin{equation}
|+, p> = {1 \over \sqrt{c_{2p} !}} A^{\dagger 2p} |+, 0>,
\quad
|-, p> = {1 \over \sqrt{c_{2p+1} !}} A^{\dagger 2p+1} |+, 0>
\end{equation}
where $|+,0>$ is the parity-even vacuum
and $c_n = n + l (1 - (-1)^n)$. 

The coherent state of $su_{q^2}(1,1)$ is defined as
\begin{equation}
|\omega; k) \equiv \sum_{p=0}^{\infty} 
{\omega^p \over \sqrt{d_p !}} |k, p>
\end{equation}
where $k=\pm$ and $d_p = [p]_{q^2}[ p + l - k{1 \over 2}]_{q^2}$.
When  $l=0$, the measure for the resolution of unity 
was given in Ref. \cite{sh}.
When $l$ is non-zero integer or half odd integer,
it was given
in terms of the $q$-deformed modified Bessel 
function  in Ref. \cite{cho}. 
For example, when $k=+1$ and $l=$ integer, we have
\begin{equation}
g(|\omega|) = {[2]_q \over 
\pi^{3 \over 2} \Pi_{j=0}^{l-1}  [j+ {1 \over 2}]_{q^2}}
|\omega|^{l - {1 \over 2}} 
\tilde K_{l - {1 \over 2}}([2]_q |\omega|)
\end{equation}
where $\tilde K_{l - {1 \over 2}}$ is the 
$q$-deformed modified Bessel function:
\begin{equation}
\tilde K_{l - {1 \over 2}}(x) = 
({x \over [2]_q})^{l - {1 \over 2}}
{ \sqrt\pi \over [n]_{q^2}!} \int_{q^{l-1}}^\infty
	d(t: q) e_q(-xt) (t^2 -1)^n_{q^2}.
\end{equation}
Once the resolution of unity is obtained for $su_{q^2}(1,1)$, 
one can easily 
demonstrate that the same holds for the oscillator 
representation.  We refer
the details to Ref. \cite{cho}.

\section{\hspace{-4mm}.\hspace{2mm} CONCLUSION}

We presented and compared various type of oscillator 
realizations of 
the $q$-deformed  $su(1,1)$ algebra and their coherent states. 
Modified Liouville measure and Bargmann 
measure are presented for the 
resolution of unity for the appropriate  coherent state. 
In addition, the suitable $q$-derivatives is  demonstrated 
for the holomorphic representation
which will be useful for the path-integral formalism 
for the $q$-deformed system.
For the application to the path integral formalism, however, 
one has to to develop the $q$-deformed higher 
dimensional integration  in terms of the non-commuting numbers 
other than the 1-dimensional one \cite{baul}.  We expect 
that this direction of research 
should accommodate the $q$-calculus 
on the hyper-plane \cite{wess}.  
Then, we may be able to apply  
the $q$-path integral formalism 
to the larger group  $SU_q(N)$ \cite{sunfu} 
and  also to the other  coadjoint orbits \cite{mpla}.

\begin{center}
{\bf ACKNOWLEDGMENTS}
\end{center}
This work is supported by the KOSEF
through the CTP at SNU and the project number
(96-1400-04-01-3, 96-0702-04-01-3),
and by the Ministry of Education through the
RIBS (BSRI/96-1419,96-2434).

\end{document}